\documentclass {article}
\usepackage[utf8]{inputenc}
\usepackage[english]{babel}
\usepackage{amsmath}
\usepackage{amsthm}
\usepackage{amssymb}

\begin{document}

    \begin{center}
        \large Alexander Zakharov \par \vspace{0.2 cm}
        \bf \Large On the rank of the intersection of free subgroups in virtually free groups
    \end{center}
    \vspace {0.5 cm}


     We prove an estimate for the rank of the intersection of free subgroups in virtually free groups, which is analogous to the Hanna Neumann inequality for subgroups in a free group and to the S.V. Ivanov estimate for subgroups in free products of groups. We also prove a more general estimate for the rank of the intersection of free subgroups in the fundamental group of a finite graph of groups with finite edge groups.

    \vspace {0.5 cm}
     \subsection* {1. Introduction.} \par
    Suppose first $G$ is a free group, $H$ and $K$ are finitely generated subgroups in $G$.
    In 1954 Howson \cite {1} proved that in this case subgroup $H \cap K$ is also finitely generated. Then
    in 1957 Hanna Neumann \cite {2} proved the following estimate for the rank of intersection of subgroups in a free group ({\it Hanna Neumann inequality}):
    \begin {equation}\label{eq1}
          \overline{r}(H \cap K) \leqslant 2 \: \overline{r}(H)\overline{r}(K),
    \end {equation}
    where $\overline{r}(H)=max \: (r(H)-1, 0)$ is the reduced rank of subgroup $H$,   $r(H)$ is the rank of subgroup $H$. \par
    In 2011 Igor Mineyev \cite{9} and Joel Friedman \cite{12} independently proved the famous Hanna Neumann conjecture which states that the coefficient 2 in the inequality (\ref{eq1}) can be omitted:
    $$\overline{r}(H \cap K) \leqslant \overline{r}(H)\overline{r}(K). $$
    \par
    S.V.Ivanov proved an estimate for subgroups of free products of groups, which is analogous to the Hanna Neumann inequality. Namely, in 1999 S.V.Ivanov \cite{3} proved that, if $G=G_{1}*G_{2}$ is a free product of groups, $H$ and $K$ are finitely generated subgroups in $G$ which intersect trivially with all the conjugates to the factors $G_{1}$ and $G_{2}$ (therefore, according to Kurosh subgroup theorem \cite{8}, $H$ and $K$ are free), then the intersection $H\cap K$ is also finitely generated and the following estimate holds:
    \begin {equation}\label{eq2}
     {\overline{r}}(H \cap K) \leqslant  6 \: {\overline{r}}(H){\overline{r}}(K).
    \end {equation}
    Later S.V.Ivanov and W.Dicks \cite{4} proved a more precise estimate for subgroups of free products which generalizes the estimate (\ref{eq2}),
    and S.V.Ivanov \cite{10} proved an analogous bound for the Kurosh rank of (arbitrary) subgroups of a free product.
     \par
    The author \cite{11} proved an estimate for the rank of the intersection of free subgroups
    in free products of groups amalgamated over a finite normal subgroup, this estimate generalizes the inequality (\ref{eq2})
    and the estimate proved by S.V.Ivanov and W.Dicks in \cite{4}.
    \par

    In this article we prove an estimate which generalizes the inequality (\ref{eq2}) to the case of subgroups of the fundamental group
    of a finite graph of groups with finite edge groups. Estimates for the rank of the intersection of subgroups in free products of groups amalgamated over a finite subgroup, as well as subgroups in HNN-extensions of groups with finite associated subgroups, follow as corollaries. Another corollary, which is obtained by applying a theorem of Stallings, is an estimate for the rank of the intersection of free subgroups in virtually free groups.\par

    \newtheorem{lemma}{Lemma}

    \subsection* {2. Bass-Serre theory.}
        In this article we use Bass-Serre theory of groups acting on trees. The main facts from this theory which we use are described below. More detailed description of this theory can be found in \cite{6}, \cite{7}. \par
       First we remind some definitions from graph theory and fix the notations. \par

        \subsubsection* {Graphs, quotient graphs} A graph $X$ is a tuple consisting of a nonempty set of vertices $V(X)$, a set of edges $E(X)$ and three mappings: $\alpha: E(X) \rightarrow V(X)$ (beginning of an edge), $\omega: E(X) \rightarrow V(X)$ (end of an edge) and $^{-1}: E(X) \rightarrow E(X)$ (inverse edge) such that $(e^{-1})^{-1}=e, \: e^{-1} \neq e, \: \alpha (e) =$ $=\omega (e^{-1})$ for every $e \in E(X)$. A graph is called finite if the sets of its edges and vertices are finite. The notion of a subgraph can be defined in a natural way. A morphism from a graph $X$ to a graph $Y$ is a map $p$ from the set of vertices and edges of $X$ to the set of vertices and edges of $Y$ which sends vertices to vertices, edges to edges and such that $\qquad \qquad \qquad \qquad \qquad \alpha(p(e))=p(\alpha(e)), \: \omega(p(e))=p(\omega(e)), \: p(e^{-1})=(p(e))^{-1}$. A bijective morphism of graphs is called an isomorphism. The degree of a vertex $v \in V(X)$ is the number of edges of the graph $X$ beginning in $v$ (notation: $deg \: v$). A morphism of graphs is called locally injective if
        it sends every two different edges beginning in the same vertex to different edges. \par
         A graph is called oriented if in each pair of its mutually inverse edges $e, e^{-1}$ one edge is chosen and called positively oriented; the other edge is called negatively oriented.
         The set of all positively oriented edges of the graph $X$ will be denoted by $E(X)^{+}$.

        A sequence $l = e_{1}e_{2}...e_{n}$ of edges of a graph $X$ is called a path beginning in $\alpha(e_{1})$ and ending in $\omega(e_{n})$ if $\omega(e_{i})=\alpha(e_{i+1}), \: i=1,...,n-1.$ (We assume that any vertex $v$ of $X$ is also a path beginning and ending in $v$, which we call trivial path at $v$.) A path is called reduced if it does not contain subpaths of the form $dd^{-1}$, where $d \in E(X)$. A path is called cyclically reduced if it is reduced and its first edge does not coincide with the inverse to its last edge; a trivial path is also cyclically reduced. A path is closed if its beginning and end coincide. A graph $X$ is called connected if for any two of its vertices $u$ and $v$ there exists a path in $X$ beginning in $u$ and ending in $v$. A tree is a connected graph which has no nontrivial reduced closed paths. A maximal subtree of a connected graph $X$ is a subtree which is maximal with respect to inclusion; it is easy to see that a maximal subtree of $X$ contains all vertices of $X$. The image of a path under a morphism of graphs is defined in a natural way: $p(e_{1}e_{2}...e_{n})=p(e_{1})p(e_{2})...p(e_{n})$.\par

        Suppose $X$ is a connected graph with a distinguished vertex $x$. Two closed paths $p_{1}$ and $p_{2}$ in $X$
         beginning in $x$ are called homotopic if  $p_{2}$ can be obtained from $p_{1}$ by a finite number of insertions and deletions of subpaths of the form $ee^{-1}$, $e \in E(X)$. One can see that the set of all reduced closed paths in $X$ beginning in $x$ forms a group with respect to the following multiplication: the product of two reduced paths $e_{1}e_{2}...e_{n}$ and $f_{1}f_{2}...f_{m}$ ($e_{i}, f_{j} \in E(X), \: i=1,...,n, \: j=1,...,m$) is the unique reduced closed path beginning in $x$ which is homotopic to the path $e_{1}e_{2}...e_{n}f_{1}f_{2}...f_{m}$; the identity of this group is the trivial path at the vertex $x$; the inverse to the path $l=e_{1}...e_{n}$ is the path $l^{-1}=e_{n}^{-1}...e_{1}^{-1}$. This group is called the fundamental group of the graph $X$ with respect to the vertex $x$ and is denoted by $\pi_{1}(X, x)$. It is easy to see that the fundamental group of a connected graph $X$ does not depend on the choice of a distinguished vertex in $X$ up to isomorphism. The isomorphism class of the fundamental group of $X$ is denoted by $\pi_{1}(X)$. \par

        The fundamental group of every connected graph $X$ is free. Moreover, let $S$ be a maximal subtree in $X$, $r_{v}$ (for each $v \in V(X)$) denote the unique reduced path beginning in $x$ and ending in $v$ which lies in $S$, $\: q_{e} = r_{\alpha(e)}er_{\omega(e)}^{-1}$ (for each $e \in E(X))$. Suppose $X$ is oriented (and thus $S$ is oriented as well). Then one can prove that the paths $q_{e}, \; e \in E(X)^{+}-E(S)^{+},$ are free generators of the group $\pi_{1}(X, x)$ (see \cite{6}).  \par
        Suppose the graph $X$ is finite. Then the following holds:
        \begin {equation}\label{new}
        r(\pi_{1}(X,x))= |E(X)^{+}|-|E(S)^{+}|=|E(X)^{+}|-|V(X)|+1,
        \end{equation}
    where the last equality holds since $S$ is a tree which contains all vertices of $X$.

        \par
       We say that a group $G$ acts on a graph $X$ on the left if left actions of $G$ on the sets $V(X)$ and $E(X)$ are defined so that $g \alpha(e)= \alpha(ge)$, $g \omega(e) = \omega(ge)$ and $ge^{-1}=(ge)^{-1}$ for all $g \in G, \: e \in E(X)$. We say that $G$ acts on $X$ without inversion of edges if $ge\neq e^{-1}$ for all $e \in E(X), \: g \in G$.  \par
       Let $G$ be a group acting on a graph $X$ without inversion of edges. For $x \in V(X) \cup E(X)$ denote by $Stab_{G} \: {x}$ the stabilizer of $x$ under the action of $G$ and by $Orb_{G} \: (x)$ the orbit of $x$ under the action of $G$: $Stab_{G} \: {x} = \{g \in G: \: gx = x \}$, \: $Orb_{G} \: (x)= \{gx, \: g\in G \}$. Define the {\it quotient graph} $\:G \: \backslash \: X$ (or $X \: / \: G$) as the graph with vertices $Orb_{G} \; (v), \; v \in V(X),$ and edges $Orb_{G} \: (e), \: e \in E(X)$; $Orb_{G} \: (v)$ is the beginning of $Orb_{G} \: (e)$ (in $G \: \backslash \: X$) if there exists $g \in G$ such that $gv$ is the beginning of $e$ (in $X$); the inverse of the edge $Orb_{G} \: (e)$ is the edge $Orb_{G} \: (e^{-1})$. \par
       Notice that the edges $Orb_{G} \: (e)$ and $Orb_{G} \: (e^{-1})$ do not coincide since $G$ acts on $X$ without inversion of edges. It is easy to see that the map $p: \: X \rightarrow G \: \backslash \: X, \quad p(x)= Orb_{G} \: (x),$ $\qquad x \in V(X) \cup E(X),$ is a surjective morphism of graphs; we call it the projection on the quotient graph. \par
       We say that $G$ acts freely on $X$ if all edge and vertex stabilizers under this action are trivial.

    \par
    \subsubsection* {The fundamental group of a graph of groups}
    A {\it graph of groups} $(\Gamma, Y)$ consists of a connected graph $Y$, vertex groups $G_{v}$ for each vertex $v \in V(Y)$, edge groups $G_{e}$ for each edge $e \in E(Y)$ such that $G_{e}=G_{e^{-1}}$ for all $e \in E(Y)$, and group embeddings
    $\alpha_{e}: G_{e} \rightarrow G_{\alpha(e)}, \; e \in E(Y)$. One can also define group embeddings $\omega_{e}: G_{e} \rightarrow G_{\omega(e)}, \: \omega_{e}=\alpha_{e^{-1}}$. A graph of groups $(\Gamma,Y)$ is called finite if the graph $Y$ is finite. A graph of groups $(\Gamma,Y)$ is called a graph of finite groups if all the vertex groups (and, therefore, all the edge groups as well) of $(\Gamma,Y)$ are finite. \par

     Let $S$ be a maximal subtree of the graph $Y$. The {\it fundamental group of the graph of groups} $(\Gamma, Y)$ with respect to the maximal subtree $S$ (notation $\pi_{1}(\Gamma,Y,S)$) is the quotient group of the free product of all vertex groups $G_{v}, \: v \in V(Y),$ and the free group with basis $\{t_{e}, \: e \in E(Y)\}$ by the normal closure of the set of the following elements:
    $$t_{e}^{-1}\alpha_{e}(g)t_{e} \cdot (\alpha_{e^{-1}}(g))^{-1} \quad(e \in E(Y), \: g \in G), \qquad t_{e}t_{e^{-1}}
    \: (e \in E(Y)), \qquad t_{e} \: (e \in E(S)). $$
    \par
    One can prove (see \cite{6}) that the fundamental group $\pi_{1}(\Gamma,Y,S)$ of the graph of groups $(\Gamma, Y)$ does not depend on the choice of the maximal subtree $S$ in $Y$ up to isomorphism. Therefore we will sometimes speak about the fundamental group of a graph of groups without mentioning the maximal subtree. One can also prove (see \cite{6}) that the vertex groups $G_{v}, \: v \in V(Y),$ can be canonically embedded in the group $\pi_{1}(\Gamma,Y,S)$.
     \par
         Consider the following examples. Suppose $Y$ consists of one pair of mutually inverse edges $e, e^{-1}$ and two vertices $u, v$ (of degree 1), then it is easy to see that the fundamental group of the graph of groups $(\Gamma, Y)$ is isomorphic to the free product of groups $G_{u}$ and $G_{v}$ amalgamated over the subgroup $\alpha_{e}(G_{e})=\omega_{e}(G_{e})$. \par
         Suppose $Y$ consists of one pair of mutually inverse edges $e, e^{-1}$ and one vertex $u$ (of degree 2), then it is easy to see that the fundamental group of the graph of groups $(\Gamma, Y)$ is isomorphic to the HNN-extension with base group $G_{u}$ and associated subgroups $\alpha_{e}(G_{e})$ and $\omega_{e}(G_{e})$. \par
         \par
          Notice that if ($\Gamma, Y$) is an arbitrary finite graph of groups and $S$ is a maximal subtree in $Y$, then the group $\pi_{1}(\Gamma, Y, S)$ can be obtained by successive applications of amalgamated free product construction (corresponding to the positively oriented edges of $S$), followed by successive applications of HNN-extension construction (corresponding to the positively oriented edges of $Y$, not belonging to $S$). \par
         One can see that if all the vertex groups (and, therefore, all the edge groups as well) of $(\Gamma, Y)$ are trivial then the fundamental group of the graph of groups ($\Gamma, Y$) is isomorphic to the fundamental group of $Y$, in particular, this group is free. Indeed, it follows from the definition of the fundamental group of a graph of groups that $\pi_{1}(\Gamma, Y, S)$ is free in this case; furthermore, its rank is equal to the number of positively oriented edges of $Y$ not belonging to the maximal subtree $S$ of $Y$, and the rank of $\pi_{1}(Y)$ is equal to the same number, as mentioned above.  \par


    \par
    \subsubsection* {Bass-Serre theorem}
    The following theorem shows the connection between fundamental groups of graphs of groups and groups acting on trees (without inversion of edges). More details and the proof of this theorem can be found in \cite{6} and \cite{7}.
    \par
    \newtheorem*{theorembs}{Theorem (Bass, Serre)}
    \begin{theorembs}
    \: $(1) \:$ Let $G=\pi_{1}(\Gamma, Y, S)$ be the fundamental group of a graph of groups $(\Gamma, Y)$ (with respect to a maximal subtree $S$). Then the group $G$ acts without inversion of edges on some tree $T$ so that
     \begin{enumerate}
     \item The quotient graph $G \: \backslash \: T$ is isomorphic to the graph $Y$.
     \item The stabilizer of a vertex $v$ under the action of $G$ is conjugate to the vertex group $G_{p(v)}$ of the graph of groups $(\Gamma, Y)$ for every $v \in V(T)$.
     \item The stabilizer of an edge $e$ under the action of $G$ is conjugate to the edge group $G_{p(e)}$ of the graph of groups $(\Gamma, Y)$ for every $e \in E(T)$.
     \end{enumerate}
     (Here $p: T \rightarrow G \: \backslash \: T$ is the projection on the quotient graph; due to the condition 1 we can identify the graphs $Y$ and $G \: \backslash \: T$ and assume that $p: T \rightarrow Y$.) \par
    $(2) \:$ Conversely, let $G$ be a group acting without inversion of edges on a tree $T$. Then the group $G$ is isomorphic to the fundamental group $\pi_{1}(\Gamma,Y,S)$ of some graph of groups $(\Gamma, Y)$ such that the conditions 1, 2, 3 from the first part of the theorem hold. In particular, each vertex group of $(\Gamma, Y)$ is equal to the stabilizer of some vertex of $T$ and each edge group of $(\Gamma, Y)$ is equal to the stabilizer of some edge of $T$.
      \end{theorembs}
      \par

    \subsubsection* {Subgroups of the fundamental group of a graph of groups which intersect trivially with the conjugates to the vertex groups}

       Suppose $G = \pi_{1}(\Gamma, Y, S)$ and $H \subseteq G$ is a subgroup which intersects trivially with the conjugates to all the vertex groups
       (and, therefore, all the edge groups as well) of $(\Gamma, Y)$. According to the first part of Bass-Serre theorem, $G$ acts on a tree        $T$ (without inversion of edges) so that conditions 1, 2, 3 of the theorem hold.\par
        Therefore, $H$ also acts on the tree $T$ in a natural way (we restrict the action of $G$ to its subgroup $H$).
         Notice that $H$ acts freely on $T$ since for each $v\in V(T)$ we have  $Stab_{H}\: (v) = Stab_{G} \: (v) \cap H = \{1\}$.
          The last equality holds since, according to the condition 2 of Bass-Serre theorem, subgroup $Stab_{G}\:(v)$ is conjugate to
          some vertex group of $(\Gamma, Y)$, and subgroup $H$ intersects trivially with all the conjugates to the vertex groups. \par
       According to the second part of Bass-Serre theorem, we obtain that $H\cong \pi_{1}(\Gamma', Y', S')$, where $(\Gamma ', Y')$ is a graph of groups, $S'$ is a maximal subtree in $Y'$, $Y'$ is isomorphic to $H \: \backslash \: T$ (due to condition 1) and all vertex groups of $(\Gamma', Y')$ are equal to the stabilizers of some vertices of $T$ under the action of $H$, and thus are trivial. Therefore, as mentioned above, $$H\cong \pi_{1}(\Gamma', Y', S') \cong \pi_{1}(Y') \cong \pi_{1}(H \: \backslash \: T),$$
      in particular, $H$ is free. \par
      Thus if subgroup $H \subseteq G=\pi_{1}(\Gamma, Y, S)$ intersects trivially with the conjugates to all the vertex groups of $(\Gamma,Y)$ then $H$ is free and, moreover,
      \begin{equation}\label{free}
      H \cong \pi_{1}(H \: \backslash \: T),
      \end{equation}
      where $T$ is a tree from Bass-Serre theorem corresponding to $G$.


    \subsection* {3. The main results.}
    \newtheorem{theorem}{Theorem}
     \begin{theorem}\label{th1}
     Suppose $G$ is the fundamental group of a finite graph of groups $(\Gamma, Y)$ with finite edge groups,
     $H, K \subseteq G$ are finitely generated subgroups which intersect trivially with the conjugates to all the vertex groups of $(\Gamma, Y)$ (and are, therefore, free). Then the following estimate holds:
     \begin{equation}\label{eq3.1}
     {\overline{r}}(H\cap K) \leqslant 6m \cdot {\overline{r}}(H)\:{\overline{r}}(K),
     \end{equation}
     where
     \begin{equation}\label{m}
      m = \max _{e \in E(Y), \: g \in G} |g^{-1}G_{e}g \cap HK|.
     \end{equation}
     In particular,
     \begin{equation}\label{eq3.11}
     {\overline{r}}(H\cap K) \leqslant 6m' \cdot {\overline{r}}(H)\:{\overline{r}}(K),
     \end{equation}
     where $m'$ is the maximum of the orders of the edge groups of $(\Gamma, Y)$.
     \end{theorem}
     (Remind that $\overline{r}(H)=max \: (r(H)-1, 0)$ is the reduced rank of subgroup $H$.) \par
     Notice that $m \leqslant m'$, so (\ref{eq3.11}) follows immediately from (\ref{eq3.1}). \par
     Applying Theorem \ref{th1} in the case when $Y$ consists of one pair of mutually inverse edges and two vertices, we obtain the following corollary.
     \newtheorem {sledstvie}{Corollary}
     \begin{sledstvie}
      Let $G$ be a free product of two groups amalgamated over a finite subgroup $T$, $H, K \subseteq G$ be finitely generated subgroups which intersect trivially with the conjugates to the factors of $G$ (and are, therefore, free). Then the following estimate holds:
     $$ {\overline{r}}(H\cap K) \leqslant 6m \cdot {\overline{r}}(H)\:{\overline{r}}(K), $$
     where $$m = \max_{g \in G} |g^{-1}Tg \cap HK|.$$
     In particular, $${\overline{r}}(H\cap K) \leqslant 6|T| \cdot {\overline{r}}(H)\:{\overline{r}}(K).$$
    \end{sledstvie}

    Applying Theorem \ref{th1} in the case when $Y$ consists of one pair of mutually inverse edges and one vertex, we obtain the following corollary.

     \begin{sledstvie}
      Let $G$ be an HNN-extension with finite associated subgroups $A_{1}, A_{2}$, and let $H, K \subseteq G$ be finitely generated subgroups which intersect trivially with the conjugates to the base group of $G$ (and are, therefore, free). Then the following estimate holds:
     $$ {\overline{r}}(H\cap K) \leqslant 6m \cdot {\overline{r}}(H)\:{\overline{r}}(K), $$
     where $$m = \max_{g \in G} |g^{-1}A_{1}g \cap HK|.$$
     In particular, $${\overline{r}}(H\cap K) \leqslant 6|A_{1}| \cdot {\overline{r}}(H)\:{\overline{r}}(K).$$

    \end{sledstvie}
    \par
    \vspace{0.2 cm}
    A group is called virtually free if it contains a free subgroup of finite index.
    Remind that a graph of finite groups is a graph of groups with finite edge and vertex groups.
    \newtheorem*{theorems}{Theorem (Stallings, \cite{5})}
     \begin{theorems}
        Suppose $G$ is a finitely generated group. Then $G$ is virtually free if and only if $G$ is the fundamental group of a finite graph of finite groups.
     \end{theorems}
    Below we show that the following theorem follows from Theorem \ref{th1} and Stallings theorem. \par

    \begin{theorem}\label{th2}
        Suppose $G$ is a virtually free group, subgroups $H, K \subseteq G$ are free and finitely generated. Then the following estimate holds:
    \begin{equation}\label{eq3.2}
     {\overline{r}}(H\cap K) \leqslant 6n \cdot {\overline{r}}(H)\:{\overline{r}}(K),
     \end{equation}
    where $n$ is the maximum of orders  $|P \cap HK|$ over all finite subgroups $P$ in $G$.

     In particular,
    \begin{equation}\label{eq3.3}
     {\overline{r}}(H\cap K) \leqslant 6n' \cdot {\overline{r}}(H)\:{\overline{r}}(K),
     \end{equation}
    where $n'$ is the minimal index of a free subgroup in $G$.
    \end{theorem}
    \subsection* {4. Proof of Theorem \ref{th2}.}
    Here we deduce Theorem \ref{th2} from Theorem \ref{th1} and Stallings theorem.
    \par
    First notice that it suffices to prove Theorem \ref{th2} for finitely generated group $G$. Indeed, instead of group $G$ we can consider group $G_{0} \subseteq G$ which is generated by subgroups $H$ and $K$; $G_{0}$ is finitely generated since $H$ and $K$ are finitely generated. Notice that $G_{0}$ is virtually free as a subgroup of a virtually free group. (Indeed, let $F \subseteq G$ be a free subgroup of finite index in $G$, then $G_{0} \cap F \subseteq G_{0}$ is a free subgroup of finite index in $G_{0}$.) It is obvious that the number $n$ from Theorem 2 will not increase when passing from $G$ to $G_{0}$. Therefore it suffices to prove the estimate (\ref{eq3.2}) for $G_{0}$.
    \par
    Thus, we can suppose that $G$ is finitely generated. Applying Stallings theorem, we obtain that $G$ is a fundamental group of a finite graph of finite groups $(\Gamma, Y)$. Moreover, subgroups $H$ and $K$ are finitely generated and intersect trivially with the conjugates to all the vertex groups of $(\Gamma, Y)$ (since $H$ and $K$ are free, and the vertex groups of $(\Gamma, Y)$ are finite). Therefore, all the conditions of Theorem \ref{th1} hold. Applying this theorem, we obtain that the estimate (\ref{eq3.1}) holds. It is obvious that the number $m$ from (\ref{m}) is less than or equal to $n$ from Theorem 2. Thus inequality (\ref{eq3.2}) holds. \par
    To prove Theorem 2 it suffices to show that $n \leqslant n'$. Moreover, the maximum of orders of finite subgroups in $G$ is less than or equal to $n'$. Indeed, otherwise there exists a finite subgroup $Q \subseteq G$ such that $|Q|>n'=|G:F|$, where the subgroup $F \subseteq G$ is free. Then there exist $g_{1}\neq g_{2} \in Q: g_{1}F=g_{2}F$, therefore $1\neq g_{2}^{-1}g_{1}\in Q \cap F$, and we get a contradiction since  $Q$ is finite, and $F$ is free. Thus estimate (\ref{eq3.3}) holds.

    \subsection* {5. Proof of Theorem \ref{th1}.}

     Applying Bass-Serre theorem, we can reformulate Theorem \ref{th1} in terms of groups acting on trees as following:
      \newtheorem*{theoremr}{Theorem $1'$}
      \begin{theoremr}
        Suppose $G$ is a group acting without inversion of edges on a tree $T$ so that the quotient graph $T/G$ is finite and all edge stabilizers  are finite. Let $H, K \subseteq G$ be finitely generated subgroups which act freely on $T$ (the action of $H$ and $K$ is restricted from the action of $G$). Then
        \begin{equation}\label{eq5.0}
     {\overline{r}}(H\cap K) \leqslant 6m \cdot {\overline{r}}(H)\:{\overline{r}}(K),
     \end{equation}
        where  $$m = \max_{x \in E(T)} |Stab_{G}(x)\cap HK|.$$

      \end{theoremr}
     We now prove Theorem $1'$. \par
     Define the projections $\pi_{H}: T/(H\cap K) \rightarrow T/H$ and $\pi_{K}: T/(H\cap K) \rightarrow T/K$ as follows:
     \begin{equation}\label{eq5.2}
        \pi_{H}(Orb_{H\cap K} x)= Orb_{H}{x}, \qquad  \pi_{K}(Orb_{H\cap K} x)= Orb_{K}{x},  \qquad x \in V(Y) \cup E(Y).
     \end{equation}
     It is easy to see that $\pi_{H}$ and $\pi_{K}$ are well-defined graph morphisms. \par
     We will now prove a few lemmas. \par

     \begin{lemma}\label{lemma 5.1}
        Suppose that the conditions of Theorem $1'$ hold.
        Then graph morphisms $\pi_{H}$ and $\pi_{K}$ are locally injective.

     \end{lemma}

    $\square$ \:    Suppose that $x_{1} = Orb_{H\cap K}(z_{1})$ and $x_{2} = Orb_{H\cap K}(z_{2})$ are two different edges of the graph $T/(H\cap K)$ beginning in a common vertex $w=Orb_{H\cap K}(u)$, where $u \in V(T)$, $z_{1}, z_{2} \in E(T)$. We can suppose that both edges $z_{1}$ and $z_{2}$ begin in $u$. Indeed, suppose that $\alpha(z_{1})=u_{1}$ and $\alpha(z_{2})=u_{2}$. Then $u_{1}=g_{1}u$ and $u_{2}=g_{2}u$, where $g_{1}, g_{2} \in H \cap K$, according to the definition of a quotient graph. Thus we can consider edges $g_{1}^{-1}z_{1}$ and $g_{2}^{-1}z_{2}$, both beginning in $u$, instead of $z_{1}$ and $z_{2}$, without changing $x_{1}$ and $x_{2}$. \par
    Now suppose that $\pi_{H}(x_{1})= \pi_{H}(x_{2})$. Then $Orb_{H}(z_{1})=Orb_{H}(z_{2})$, or $z_{1}=hz_{2}$, where $h \in H$. Therefore, $u=hu$, but $H$ acts freely on $T$, so $h = 1$. Thus, $z_{1}=z_{2}$ and $x_{1}=x_{2}$, and we obtain a contradiction. This shows that $\pi_{H}$ is locally injective. Analogously $\pi_{K}$ is locally injective.

           $\blacksquare$
     \begin{lemma}\label{lemma 5.2a}
        Suppose a group $G_{0}$ acts on a set $M$, $H_{0}, K_{0}$ are subgroups of $G_{0}$, and $z \in M$. Then $Orb_{H_{0}}(z) \, \cap \, Orb_{K_{0}}(z)$ consists of not more than $|Stab_{G_{0}}(z)\cap H_{0}K_{0}|$ orbits under the action of $H_{0} \cap K_{0}$. \par
     \end{lemma}

     $\square$ \: 
     Suppose that $u_{1}=h_{1}z=k_{1}z$, \: $u_{2}=h_{2}z=k_{2}z$, ..., \: $u_{n}=h_{n}z=k_{n}z$ are $n$ distinct elements of $Orb_{H_{0}}(z) \, \cap \, Orb_{K_{0}}(z)$, where $n > |Stab_{G_{0}}(z)\cap H_{0}K_{0}|$ and $h_{i} \in H_{0}, \: k_{i} \in K_{0}, \: i = 1,2,...,n$. It suffices to show that at least two of $u_{1},u_{2},...,u_{n}$ belong to the same orbit under the action of $H_{0} \cap K_{0}$. \par
     Indeed, $h_{1}^{-1}k_{1}, \: h_{2}^{-1}k_{2},..., h_{n}^{-1}k_{n} \in Stab_{G_{0}}(z)\cap H_{0}K_{0}$, thus there exist $p, q \in \{1, 2,..., n\}$ such that $h_{p}^{-1}k_{p}=h_{q}^{-1}k_{q}$. Then $h_{q}h_{p}^{-1}=k_{q}k_{p}^{-1}=c \in H_{0} \cap K_{0}$ and $u_{q}=h_{q}z=ch_{p}z=cu_{p}$. Thus $u_{p}$ and $u_{q}$ belong to the same orbit under the action of $H_{0} \cap K_{0}$, and Lemma \ref{lemma 5.2a} is proven.

     $\blacksquare$

     \begin{lemma}\label{lemma 5.2}
        Suppose that the conditions of Theorem $1'$ hold. Let $e, f$ be edges of the graphs $T/H$, $T/K$ respectively.
         Then the number $N$ of edges of the graph $T/(H\cap K)$ which project under $\pi_{H}$ into $e$ and under $\pi_{K}$ into $f$ (simultaneously) is not bigger than $m$, where
        $m = max \: (|Stab_{G}(x)\cap HK|, \: x \in E(T))$.
     \end{lemma}

       $\square$ \: Suppose $e = Orb_{H}(r), \: f = Orb_{K}(s), \: r, s \in E(T)$. \par
       Notice that $N$ is equal to the number of orbits $Orb_{H\cap K}(t), \: t \in E(T),$ such that $\pi_{H}(Orb_{H\cap K}(t))=Orb_{H}(r)$ and $\pi_{K}(Orb_{H\cap K}(t))=Orb_{K}(s)$, or, equivalently, $t \in Orb_{H}(r) \cap Orb_{K}(s)$. If $Orb_{H}(r) \cap Orb_{K}(s) = \varnothing$, then Lemma \ref{lemma 5.2} holds. Otherwise, let $d \in Orb_{H}(r) \cap Orb_{K}(s), \: d \in E(T).$ Then $Orb_{H}(r)=Orb_{H}(d)$ and $Orb_{K}(s) = Orb_{K}(d)$.
       Apply Lemma \ref{lemma 5.2a} with $G_{0}=G, \: H_{0}=H, \: K_{0}=K, \: M=E(T), \: z = d$. This proves Lemma \ref{lemma 5.2}.

    $\blacksquare$

         \par \vspace{0.2 cm}
     Notice that if the conditions of Theorem $1'$ hold and the graph $T/H$ is a tree, then, since $H \cong \pi_{1}(T/H)$ according to (\ref{free}), the subgroup $H$ is trivial and therefore $H\cap K$ is trivial as well. In this case the estimate (\ref{eq5.0}), and thus Theorem $1'$, hold; the same is true if the graph $T/K$ or $T/(H\cap K)$ is a tree. Thus we can assume below that the graphs $T/H$, $T/K$ and $T/(H\cap K)$ are not trees. \par

     Suppose a connected graph $X$ is not a tree. We call the {\it core} of the graph $X$ a subgraph of $X$ which consists of all vertices and edges of $X$ which belong to any nontrivial closed cyclically reduced path in $X$. Notice that if a vertex $v \in V(X)$ belongs to the core of $X$ then the core of $X$ consists of all vertices and edges of $X$ which belong to any reduced closed path in $X$ beginning in $v$. \par
      For a nontrivial subgroup $H \subseteq G$ denote by $\Psi(H)$ the core of the graph $T/H$. Notice that since $T/H$ is connected $\Psi(H)$ is also connected. Notice also that the graph $\Psi(H)$ does not contain vertices of degree less than 2. We fix an arbitrary orientation of $\Psi(H)$. \par

    \begin{lemma}\label{lemma 5.3}
        Suppose that the conditions of Theorem $1'$ hold (in particular, subgroup $H \subseteq G$ is finitely generated and acts freely on $T$) and subgroup $H$ is nontrivial. Then the graph $\Psi(H)$ is finite, $H \cong \pi_{1}(\Psi(H))$ and the following equalities hold:
        \begin{equation}\label{eq5.3}
         \overline{r}(H) = |E(\Psi(H))^{+}|-|V(\Psi(H))|=\frac{1}{2}\sum_{v \in V(\Psi(H))} \: (deg \: v - 2).
        \end{equation}
        Similar statement holds for subgroups $K$ and $H\cap K$ from Theorem $1'$.
    \end{lemma}

    $\square$ \:
    It was shown above (see (\ref{free})) that $H\cong \pi_{1}(T/H)$ since $H$ acts freely on $T$ and due to Bass-Serre theorem. Fix a vertex $v \in V(T/H)$ which lies in the subgraph $\Psi(H)$. Since any reduced closed path in $T/H$ beginning in $v$ lies in $\Psi(H)$, we obtain that $\pi_{1}(T/H, v) \cong \pi_{1}(\Psi(H), v)$, therefore, $H \cong \pi_{1}(\Psi(H))$. \par

    Suppose reduced paths $p_{1},...,p_{n}$ are free generators of the group $\pi_{1}(\Psi(H), v)$; $n$ is finite, since $H \cong \pi_{1}(\Psi(H))$ and $H$ is finitely generated. Any reduced closed path in $\Psi(H)$ beginning in $v$ is a product of some paths from $p_{1},...,p_{n}$ and their inverses. As mentioned above, any edge $e$ of the graph $\Psi(H)$ belongs to some closed reduced path in $\Psi(H)$ beginning in $v$, so $e$ belongs to at least one of the paths $p_{1},...,p_{n}$ and their inverses. Thus, the graph $\Psi(H)$ is finite. \par

    Moreover, according to (\ref{new}), we obtain:

    $$r(H)=|E(\Psi(H))^{+}|-|V(\Psi(H))|+1.$$
    Therefore, the first equality in (\ref{eq5.3}) holds. \par
    Finally, the sum of the degrees of all vertices of any (oriented) graph is equal to the doubled number of its positively oriented edges,  therefore, the second equality in (\ref{eq5.3}) holds as well. \par
    It is obvious that the same proof is true for the subgroups $K$ and $H\cap K$ from Theorem $1'$. $\blacksquare$

    \begin{lemma}\label{lemma 5.4}
        Suppose that the conditions of Theorem $1'$ hold and $H\cap K$ is nontrivial. Then the image of the graph $\Psi(H\cap K)$ under the projection $\pi_{H}, \pi_{K}$ lies in the graph $\Psi(H)$, $\Psi(K)$ respectively. Thus, we can consider the restriction of the projections $\pi_{H}: \Psi(H\cap K) \rightarrow \Psi(H)$, \qquad \qquad \qquad \qquad \qquad $\pi_{K}: \Psi(H\cap K) \rightarrow \Psi(K)$.
    \end{lemma}

    $\square$ \:
    This lemma follows from Lemma \ref{lemma 5.1}. Indeed, suppose $v$ is a vertex of the graph $\Psi(H\cap K)$. Then $v$ belongs to some closed cyclically reduced path $p$ in $T/(H\cap K)$. Since $\pi_{H}$ is locally injective (due to Lemma \ref{lemma 5.1}), the closed path $\pi_{H}(p)$ in $T/H$ is also cyclically reduced, and $\pi_{H}(v)$ belongs to this path, therefore, $\pi_{H}(v)$ belongs to the graph $\Psi(H)$. Similarly $\pi_{K}(v)$ belongs to the graph $\Psi(K)$. The same is true for edges of $\Psi(H \cap K)$. \:     $\blacksquare$

    \begin{lemma}\label{lemma 5.5}
        Suppose that the conditions of Theorem $1'$ hold and $H\cap K$ is nontrivial. Let $a, b$ be vertices of the graphs $\Psi(H)$, $\Psi(K)$ respectively. Let $w_{1},...,w_{s}$ be all vertices of the graph $\Psi(H\cap K)$ which project under $\pi_{H}$ into $a$ and under $\pi_{K}$ into $b$. Then the following inequalities hold:

        \begin{equation}\label{eq5.4}
            deg \: w_{i} \leqslant deg \: a, \quad deg \: w_{i} \leqslant deg \: b, \quad i=1,...,s
        \end{equation}

        \begin{equation}\label{eq5.5}
            \sum_{i=1}^{s} \: deg \: w_{i} \leqslant m \cdot deg \: a \cdot deg \: b,
        \end{equation}
        where $m = max \: (|Stab_{G}(x)\cap HK|, \: x \in E(T)).$
        (We define $s$ and the sum in (\ref{eq5.5}) to be 0 if there are no such vertices $w_{i}$.)
    \end{lemma}
    $\square$ \:
    Due to Lemma \ref{lemma 5.4} each edge of the graph $\Psi(H\cap K)$ beginning in one of the vertices
    $w_{i} \: (i = 1, ... , s)$ projects under $\pi_{H}$ into an edge of the graph $\Psi(H)$ beginning in $a$, and projects under $\pi_{K}$ into an edge of the graph $\Psi(K)$ beginning in $b$. \par
    Inequality (\ref{eq5.4}) now follows immediately from Lemma \ref{lemma 5.1}. \par
    Applying Lemma \ref{lemma 5.2} for every edge $x$ of the graph $\Psi(H)$ beginning in
    $a$ and for every edge $y$ of the graph $\Psi(K)$ beginning in $b$, we obtain inequality (\ref{eq5.5}) as well. \: $\blacksquare$

    \par
    \vspace{0.3 cm}

    We will now complete the proof of Theorem $1'$. The following part of the proof follows the idea of S.V.Ivanov \cite{3}. \par

    Applying the equalities (\ref{eq5.3}) from Lemma \ref{lemma 5.3}, we can reformulate the estimate (\ref{eq5.0}) of Theorem $1'$ in terms of the degrees of vertices of the graphs $\Psi$:

    \begin {equation}\label{eq11}
     \sum _{w \in V(\Psi(H\cap K))}(deg \: w - 2) \leqslant 3 m \cdot \! \! \! \! \! \! \! \sum _{a \in V(\Psi(H))}(deg \: a - 2) \cdot \! \! \! \! \! \!  \sum _{b \in V(\Psi(K))}(deg \: b - 2)
    \end {equation}

    Thus, to prove Theorem $1'$ it suffices to prove the inequality (\ref{eq11}).
    Notice that to prove the inequality (\ref{eq11}) it suffices to prove the following inequality:

    \begin {equation}\label{eq12}
     \sum _{i=1}^{s_{a,b}}(deg \: w^{a,b}_{i} - 2) \leqslant 3m \cdot (deg \: a - 2) \cdot (deg \: b - 2),
    \end {equation}
    for all vertices $a \in V(\Psi(H))$, $b \in V(\Psi(K))$. Here $w^{a,b}_{1},...,w^{a,b}_{s_{a,b}}$ are all vertices of the graph $\Psi(H\cap K)$ which project under $\pi_{H}$ into $a$ and under $\pi_{K}$ into $b$. (We define $s_{a,b}$ and the sum in (\ref{eq12}) to be 0 if there are no such vertices $w_{i}^{a,b}$.)  \par
    Indeed, suppose the inequality (\ref{eq12}) holds. Then, according to Lemma \ref{lemma 5.4}, we obtain:
    $$ \sum_{w \in V(\Psi(H\cap K))}(deg \: w - 2) = \sum_{(a,b)} \sum _{i=1}^{s_{a,b}}(deg \: w^{a,b}_{i} - 2) \leqslant \sum_{(a,b)} 3m \cdot (deg \: a - 2) \cdot (deg \: b - 2) = $$ $$= 3 m \cdot \! \! \! \! \! \! \! \sum _{a \in V(\Psi(H))}(deg \: a - 2) \cdot \! \! \! \! \! \!  \sum _{b \in V(\Psi(K))}(deg \: b - 2), $$
    where the sum $\sum_{(a,b)}$ extends over all vertices $a \in V(\Psi(H))$, $b \in V(\Psi(K))$.
     Thus, if (\ref{eq12}) holds, then (\ref{eq11}) holds as well.
         \par
    It suffices to prove the inequality (\ref{eq12}). We can assume without loss of generality that
    \begin {equation}\label{eq13}
    deg \: a \leqslant deg \: b.
    \end {equation}
    Applying Lemma \ref{lemma 5.5}, we get the following inequalities:

    \begin {equation}\label{eq14}
       deg \: w^{a,b}_{i} \leqslant deg \: a, \quad i=1,...,s_{a,b},
    \end {equation}
    \begin {equation}\label{eq15}
       \sum_{i=1}^{s_{a,b}} deg \: w^{a,b}_{i} \leqslant m \cdot deg \: a \cdot deg \: b.
    \end {equation}

    Consider two cases. If $s_{a,b} \leqslant m \cdot deg \: b$, then, applying inequality (\ref{eq14}), we obtain:
    $$ \sum _{i=1}^{s_{a,b}}(deg \: w^{a,b}_{i} - 2) \leqslant s_{a,b} (deg \: a - 2) \leqslant m \cdot deg \: b \cdot (deg \: a - 2).$$
    If $s_{a,b} \geqslant m \cdot deg \: b$, then, applying inequality (\ref{eq15}), we obtain:
    $$ \sum _{i=1}^{s_{a,b}}(deg \: w^{a,b}_{i} - 2) = \sum _{i=1}^{s_{a,b}} deg \: w^{a,b}_{i} - 2s_{a,b} \leqslant m \cdot deg \: a \cdot deg \: b - 2 m \cdot deg \: b = m \cdot deg \: b \cdot (deg \: a - 2). $$
    Thus, in any case the following inequality holds:
    \begin{equation}\label{eq17}
     \sum _{i=1}^{s_{a,b}}(deg \: w^{a,b}_{i} - 2) \leqslant m \cdot deg \: b \cdot (deg \: a - 2).
    \end{equation}
     Moreover, the graph $\Psi(H)$ has no vertices of degree less than 2. Therefore, according to (\ref{eq13}), $deg \: b \geqslant deg \: a \geqslant 2$. If $deg \: b = 2$, then $deg \: a = 2$, so (\ref{eq17}) implies (\ref{eq12}). Otherwise, if $deg \: b \geqslant 3$, then $deg \: b \leqslant 3(deg \: b - 2)$, therefore, $$m \cdot deg \: b \cdot (deg \: a - 2) \leqslant 3 m \cdot (deg \: a - 2) \cdot (deg \: b - 2),$$
    so (\ref{eq17}) again implies (\ref{eq12}). Therefore, in any case the inequality (\ref{eq12}) holds. \par
    This shows that Theorem $1'$ holds. \par
    Thus, Theorem \ref{th1} is proven.
    \par
    \vspace{0.5 cm}
    The author thanks A.A. Klyachko for many useful conversations and many useful remarks. The author also thanks W. Dicks for useful remarks.

     \vspace{0.3 cm}
    Faculty of Mechanics and Mathematics, Moscow State University, Russia \par
    {\it E-mail address:} zakhar.sasha@gmail.com \par \par
    This work was supported by the Russian Foundation for Basic Research, project no. 11-01-00945.

\end{document}